\documentclass[a4paper,12pt,reqno]{amsart}
\usepackage{graphics}
\newtheorem{lem}{Lemma}

\newtheorem{theo}{Theorem}
\newtheorem{rem}{Remark}
\newtheorem{cor}{Corollary}

\numberwithin{equation}{section}

\newcommand{\sgn}{\operatorname{sgn}}
\newcommand{\norm}{\operatorname{norm}}

\begin{document}

\title[Polynomial enumeration formulas]{A method for proving polynomial enumeration formulas}
\author[Ilse Fischer]{\box\Adr}

\newbox\Adr
\setbox\Adr\vbox{ \centerline{ \large Ilse Fischer} \vspace{0.3cm}
\centerline{Institut f\"ur Mathematik, Universit\"at Klagenfurt,}
\centerline{Universit\"atsstrasse 65 -- 67, A-9020 Klagenfurt, Austria.}
\centerline{E-mail: {\tt Ilse.Fischer@uni-klu.ac.at}} }

\date{}
\maketitle

\begin{abstract}
We present an elementary method for proving enumeration
formulas which are polynomials in certain parameters if others are fixed and
factorize into distinct linear factors over $\mathbb{Z}$. Roughly speaking the
idea is to prove such formulas by ``explaining'' their zeros using an
appropriate combinatorial extension of the objects under consideration to negative
integer parameters. We apply this method to prove a new refinement of the
Bender-Knuth (ex-)Conjecture, which easily implies the Bender-Knuth
(ex-)Conjecture itself. This is probably the most elementary way to prove
this result currently known. Furthermore we adapt our method to
$q$-polynomials, which allows us to derive generating
function results as well. Finally we use this method to give another proof for the enumeration
of semistandard tableaux of a fixed shape, which is opposed to the
Bender-Knuth (ex-)Conjecture refinement a multivariate application of our method. 
\end{abstract}

\section{Introduction}

\subsection{A simple example}
Let $F(r,k)$ denote the number of partitions 
$(\lambda_{1},\ldots,\lambda_{r})$, i.e. $\lambda_{1} \ge \lambda_{2} \ge \ldots
\ge \lambda_{r}$, of length $r$, with parts in $\{0,1,\ldots,k\}$. It
is basic combinatorial knowledge that 
$$
F(r,k)=\binom{k+r}{r}=\frac{(k+1) \cdot (k+2)\cdot \ldots \cdot (k+r)}{r!}.
$$
For fixed $r$ this expression is a polynomial in $k$ with distinct
integer zeros. In this paper we present an elementary method for proving
polyomial enumeration formulas of that type, together with some
non-trivial applications. The underlying idea is to find
the appropriate extension of the combinatorial objects under
consideration to (typically) negative integer parameters 
and with this ``explain'' the zeros of the enumeration polynomial.  

To be more concrete let us first demonstrate this $3$-step-method on the basis of our
simple example.

(1) In the first step we extend the combinatorial interpretation of $F(r,k)$
  to negative integer $k$'s. For $k < 0$ we define
$$
F(r,k) = (-1)^{r}  [ \# (\lambda_{1},\ldots,\lambda_{r}) \in \mathbb{Z}^{r} \text{ with } k < \lambda_{1} <
\lambda_{2} < \ldots < \lambda_{r} < 0].
$$
This definition seems to appear from nowwhere,
however, the following step should convince us that it was a wise choice.

(2) In this step we show that for fixed $r$ the function $k \to F(r,k)$ can be expressed
  by a polynomial in $k$ of degree at most $r$. This is equivalent to
  $\Delta^{r+1} F(r,k) =0$, where the differences are taken with respect to
  the parameter $k$. In order to show this we use induction with respect to
  $r$. The initial step follows from $F(1,k)=k+1$. Assume that $r > 1$ and $k
  \ge 0$. Then 
\begin{multline*}
\Delta F(r,k) = F(r,k+1) - F(r,k)   \\
= [ \# (\lambda_{1}, \lambda_{2}, \ldots, \lambda_{r}) \text{ with } k+1 \ge \lambda_{1} \ge
\ldots \ge \lambda_{r} \ge 0 ] \\ - [ \# (\lambda_{1},\lambda_{2},\ldots,\lambda_{r}) \text{ with } k \ge \lambda_{1} \ge
\ldots \ge \lambda_{r} \ge 0 ] \\
= [\# (\lambda_{1},\lambda_{2},\ldots,\lambda_{r}) \text{ with } k+1 \ge \lambda_{1} \ge
\ldots \ge  \lambda_{r} \ge 0 \text{ and } \lambda_{1}=k+1 ] \\
= F(r-1,k+1)
\end{multline*}
If $k < 0$ we have 
\begin{multline*}
\Delta F(r,k) = F(r,k+1) - F(r,k)   \\
= (-1)^{r} [\# (\lambda_{1},\lambda_{2},\ldots,\lambda_{r}) \text{ with } k+1 < \lambda_{1} <
\ldots < \lambda_{r} < 0] \\ - (-1)^{r} [\# (\lambda_{1},\lambda_{2},\ldots,\lambda_{r}) \text{ with } k < \lambda_{1} <
\ldots < \lambda_{r} < 0] \\
= (-1)^{r-1} [\# (\lambda_{1},\lambda_{2},\ldots,\lambda_{r}) \text{ with } k < \lambda_{1} <
\ldots < \lambda_{r} < 0 \text{ and } \lambda_{1}=k+1] \\
= F(r-1,k+1)
\end{multline*}
The induction hypothesis implies $\Delta^{r} F(r-1,k+1)=0$ and thus
$\Delta^{r+1} F(r,k)=0$.

(3) In the final step we explore the integer zeros of $F(r,k)$ in
  $k$. Consider the definition of $F(r,k)$ for negative $k$'s and observe that 
$F(r,k)=0$ for $k=-1,-2,\ldots,-r$. By Step~2 $F(r,k)$ is a polynomial in $k$ and
therefore it has the factor $(k+1)_{r}$, where the Pochhammer symbol $(a)_{n}$
is defined by $(a)_{n} = \prod_{i=0}^{n-1} (a+i)$. The degree estimation of
Step~2 implies that this factor determines $F(r,k)$ up to a factor independent
of $k$. Observe that $F(r,0)=1$, and  thus this factor is equal to $1/r!$ and
the formula is proved.

\subsection{The method} We summarize the  general strategie in
the example above and with this establish our method for proving polynomial
enumeration formulas. It applies to the enumeration of combinatorial objects which depend on an
integer parameter $k$ and where we suspect the existence of an enumeration
formula which is polynomial in $k$ and factorizes into
distinct linear factors over $\mathbb{Z}$. The method is divided into the
following three steps.

\begin{enumerate}

\item {\bf Extension of the combinatorial interpretation.}
Typically the admissible domain of $k$ is a set $S$ of non-negative integers. In
the first step of our method we have to find (most likely new) combinatorial
objects indexed by an {\it arbitrary} integer $k$ which are in
bijection with the original objects for $k \in S$.

\smallskip

\item {\bf The extending objects are enumerated by a polynomial.}
The extension of the combinatorial interpretation in the previous
step has to be chosen so that we are able to prove that the new
objects are enumerated by a polynomial in $k$. In many cases this is done with
the help of a recursion. Moreover the degree of this polynomial has to be computed. 

\smallskip

\item {\bf Exploring ``natural'' linear factors.} 
Finally one has to find the $k$'s for which there exist none of
these objects, i.e. one has to compute the (integer) zeros of the
polynomial.\footnote{In the first step it may have been necessary to introduce a 
signed enumeration outside of the admissible domain in order to have the same 
enumeration polynomial for all $k$'s. In this case we have to find
the $k$'s for which objects cancel in pairs with respect to
the sign.} Typically these zeros will not lie in $S$, which made
the extension in Step~1 necessary. Moreover one has to find a
non-zero evaluation of the polynomial which is easy to compute,
and together with the zeros the polynomial is finally computable.

\end{enumerate}

The last step shows the limits of this method. Even if one succeeds
in the first two steps, it may be that the polynomial has
non-integer zeros or multiple zeros and the method as described
does not work. On the other hand the enumeration problems
which result in polynomials that factorize totally over
$\mathbb{Z}$ are exactly the one we are especially interested in and where
we are longing for an understanding of the simplicity of the
result.

\subsection{A refinement of the Bender-Knuth (ex-)Conjecture.} Next we explain
a plane partition enumeration result we have obtained by using this method. The
main purpose of the rest of the paper is the proof of this result.
Let $\lambda=(\lambda_{1},\lambda_{2},\dots,\lambda_{r})$ be a
partition. A {\it strict plane partition} of shape $\lambda$ is an
array $\pi_{1 \le i \le r, 1 \le j \le \lambda_i}$ of non-negative
integers such that the rows are weakly decreasing and the columns
are strictly decreasing. 
The {\it norm} $n(\pi)$ of a strict
plane partition is defined as the sum of its parts and $\pi$ is
said to be a strict plane partition of the non-negative integer
$n(\pi)$. For instance

\begin{center}
\begin{tabular}{cccccc}
$7$ & $5$ & $5$ & $4$ & $3$ & $2$ \\
$6$ & $4$ & $3$ & $2$ &   &   \\
$5$ & $2$ &   &   &   &   \\
$3$ & $1$ &   &   &   &
\end{tabular}
\end{center}

{\noindent is a strict plane partition of shape $(6,4,2,2)$ with norm $52$. 
In \cite[p.50]{benderknuth} Bender and Knuth had conjectured
that the generating function of strict
plane partitions with at most $c$ columns, parts in
$\{1,2,\dots,n\}$ and with respect to this norm is equal to
$$
\sum q^{n(\pi)} = \prod_{i=1}^n \frac{[c+i;q]_i}{[i;q]_i},
$$
where $[n;q] = 1 + q + \dots + q^{n-1}$ and $[a;q]_n =
\prod_{i=0}^ {n-1} [a+i;q]$. This conjecture was proved by
Andrews~\cite{andrews}, Gordon~\cite{gordon}, Macdonald~\cite[Ex.
19, p.53]{macdonald} and Proctor~\cite[Prop. 7.2]{proctor}. For
related papers, which mostly include generalizations of the
Bender-Knuth (ex-)Conjecture see~\cite{des1, des2, kadell, kratt,
proctor2, stem2}.

Using a ``$q$-extension'' of our method we have obtained the following new
refinement of this result. As an additional parameter $k$ we introduce the number of parts
equal to $n$ in the strict plane partition.

\begin{theo}
\label{main} The generating function of strict plane partitions
with parts in $\{1,2,\dots,n\}$, at most $c$ columns and $k$
parts equal to $n$ is
$$
\sum q^{n(\pi)} =
\frac{q^{k n} [k+1;q]_{n-1} [1+c-k;q]_{n-1}}{[1;q]_{n-1}}
\prod_{i=1}^{n-1} \frac{[c+i+1;q]_{i-1}}{[i;q]_{i}}.
$$
\end{theo}

If we sum this generating function over all $k$'s, $0 \le k \le
c$, we easily obtain the Bender-Knuth (ex-)Conjecture. Probably
this detour via Theorem~\ref{main} is the easiest and most
elementary way to prove the Bender-Knuth (ex-)Conjecture currently known.
In \cite[Sec. 3]{kratt2} the authors come to the conclusion that all other
proofs of the Bender-Knuth (ex-)Conjecture ``share more or less
explicitly an identity, which relates Schur functions and odd
orthogonal characters of the symmetric group of rectangular
shape''. In our elementary proof this is not the case.

In order to illustrate our method we first prove the special case $q=1$ of
Theorem~\ref{main}, i.e. we compute the number of strict plane
partitions with parts in $\{1,2,\dots,n\}$, at most $c$ columns
and $k$ parts equal to $n$, see Theorem~\ref{special}. (Observe that for $q=1$ the formula
in Theorem~\ref{main} is a polynomial in $k$, which factorizes
into distinct linear factors over $\mathbb{Z}$.) This result is new as
well. Later we will see that the method can be extended to $q$-polynomials 
in order to prove the general result.

\subsection{Outlook and outline of the paper}
We plan to apply this method to other enumeration problems in the
future. The most ambitious project in this direction is probably
our current effort to give another proof of the refined alternating
sign matrix Theorem. There is some hope for a proof along the lines
of the proof of Theorem~\ref{main}. Let $A(n,k)$ denote the number of alternating sign
matrices of order $n$, where
the unique $1$ in the first row is in the $k$-th column. It
came as a surprise that  the number of strict
plane partitions with parts in $\{1,2,\dots,n\}$, at most $n-1$
columns and $k-1$ parts equal to $n$ divided by $A(n,k)$ is independent of
$k$. In other words: The enumeration polynomial is -- up to a constant and 
up to a shift -- equal to the enumeration polynomial in Theorem~\ref{main} if we 
set $q=1$ and $c=n-1$ there.  Thus
an application of our method to alternating sign matrices could be
very similar to the application to strict plane partitions which
is under consideration in this paper, see Section~\ref{final}.
Moreover we plan to extend our method to polynomial enumeration formulas 
that do not factor into distinct linear factors over $\mathbb{Z}$. For instance 
polynomial enumeration formulas that are certain sums of polynomials that
factorize into distinct linear factors over $\mathbb{Z}$ could be a first goal.

\medskip

The rest of the paper is organized as follows. In
Section~\ref{pattern} we introduce a combinatorial extension with respect
to $k$ of
strict plane partitions with parts in $\{1,2,\dots,n\}$, at most
$c$ columns and $k$ parts equal to $n$ as proposed in {\it Step
1} of our method above. In Section~\ref{polynomial} we show that these
objects are enumerated by a polynomial in $k$ which is of degree
$2n-2$ at most ({\it Step 2}) and in Section~\ref{zeros} we show that the 
polynomial has the predicted zeros ({\it Step 3}). This concludes the proof of
Theorem~\ref{main} for $q=1$. In Section~\ref{semistand} we apply
the method to give another proof of the formula for the number of
semistandard tableaux of a fixed shape. This application of our
method is of interest since in this case we have to work with more
than just one polynomial parameter. Finally we extend our method
to what we call ``$q$-polynomials'' and prove Theorem~\ref{main} in
its full strength in Section~\ref{q-p}. In Section~\ref{final} a
connection of our result to the refined alternating sign matrix
Theorem is presented.

\medskip

Throughout the whole article we use the extended definition of
the summation symbol, namely,
\begin{equation}
\label{sumext}
\sum_{i=a}^{b} f(i) =
\begin{cases}
f(a) + f(a+1) + \dots + f(b) & \text{if $a \le b$} \\
0 & \text{if $b=a-1$} \\
- f(b+1) - f(b+2) - \dots - f(a-1) & \text{if $b+1 \le a-1 $}
\end{cases}.
\end{equation}
This assures that for any polynomial $p(X)$ over an arbitrary integral domain $I$
containing $\mathbb{Q}$ there exists a unique polynomial $q(X)$ over $I$ such that $\sum_{x=0}^y p(x) = q(y)$ for all
integers $y$. We usually write $\sum_{x=0}^ y p(x)$ for $q(y)$.

\section{From strict plane partitions to generalized $(n-1,n,c)$ Gelfand-Tsetlin-patterns}
\label{pattern}

Let $n$, $c$ be integers, $n$ positive and $c$ non-negative. A
Gelfand-Tsetlin-pattern with $n$ rows is 
a triangular array of integers, say  
\begin{center}
\begin{tabular}{ccccccccccccccccc}
  &   &   &   &   &   &   &   & $a_{n,n}$ &   &   &   &   &   &   &   & \\ 
  &   &   &   &   &   &   & $a_{n-1,n-1}$ &   & $a_{n-1,n}$ &   &   &   &   &   &   & \\
  &   &   &   &   &   & $\dots$ &   & $\dots$ &   & $\dots$ &   &   &   &   &   & \\
  &   &   &   &   & $a_{3,3}$ &   & $\dots$ &   & $\dots$ &   & $a_{3,n}$ &   &   &   &   & \\
  &   &   &   & $a_{2,2}$ &   & $a_{2,3}$ &  &   $\dots$ &   & $\dots$   &  & $a_{2,n}$  &   &   &   & \\
  &   &   & $a_{1,1}$ &   & $a_{1,2}$ &   & $a_{1,3}$ &   & $\dots$ &   & $\dots$ &   & $a_{1,n}$ &   &   & 
\end{tabular},
\end{center}
such that $a_{i,j} \le a_{i-1,j}$ for $1 < i \le j \le n$ and 
$a_{i,j} \le a_{i+1,j+1}$ for $1 \le i \le j < n$, see~\cite[p. 313]{stan} or
\cite[(3)]{gel} for the original reference. An example of a
Gelfand-Tsetlin-pattern with $7$ rows is given below.
\begin{center}
\begin{tabular}{ccccccccccccccccc}
  &   &   &   &   &   &   &   & 1 &   &   &   &   &   &   &   & \\ 
  &   &   &   &   &   &   & 1 &   & 1 &   &   &   &   &   &   & \\
  &   &   &   &   &   & 1 &   & 1 &   & 3 &   &   &   &   &   & \\
  &   &   &   &   & 0 &   & 1 &   & 2 &   & 4 &   &   &   &   & \\
  &   &   &   & 0 &   & 1 &   & 1 &   & 3 &   & 5 &   &   &   & \\
  &   &   & 0 &   & 0 &   & 1 &   & 2 &   & 4 &   & 6 &   &   & \\ 
  &   & 0 &   & 0 &   & 0 &   & 2 &   & 2 &   & 4 &   & 6 &   & 
\end{tabular} 
\end{center}
The following correspondence between Gelfand-Tsetlin-patterns and strict plane partitions is crucial 
for our paper.

\begin{lem}
\label{bijection}
There is a bijection between Gelfand-Tsetlin-patterns $(a_{i,j})$ with $n$
rows, parts in $\{0,1,\ldots,c\}$ and fixed $a_{n,n}=k$, and strict plane partitions with parts in
$\{1,2,\dots,n\}$, at most $c$ columns and $k$ parts equal to $n$.
\end{lem}

{\it Proof.}
Given such a Gelfand-Tsetlin-pattern, the corresponding strict plane
partition is such that the shape filled by entries greater than
$i$ corresponds to the partition given by the $(n-i)$-th row of
the Gelfand-Tsetlin-pattern, the top row being the first row. As an example consider the strict plane
partition in the introduction. If we
choose $n=7$ and $c=6$ then this strict plane partition corresponds to the
Gelfand-Tsetlin pattern above. \qed

\medskip

Therefore it suffices to enumerate Gelfand-Tsetlin-patterns
$(a_{i,j})$ with $n$ rows, parts in $\{0,1,\ldots,c\}$ and fixed 
$a_{n,n}=k$. Why should this be easier than enumerating the 
corresponding strict plane partitions?

Recall that $k$ is the polynomial parameter in our refinement of the 
Bender-Knuth (ex-~)Conjecture we want to make use of when applying our 
method. In order to accomplish Step~1 of the method we have to find a ``natural'' extended
definition of  strict plane partitions with 
parts in $\{1,2,\ldots,n\}$, at most $c$ columns and $k$ parts 
equal to $n$, where $k$ is an {\it arbitrary} integer which does not necessarily
lie in $\{0,1,\ldots,c\}$. (A priori parts equal to $n$ may only appear in the first row of the 
(column-)strict plane partition with parts in $\{1,2,\ldots,n\}$ and thus $k
\in \{0,1,\ldots,c\}$.)  ``Natural'' stands for the
fact that the extension has to be chosen such that the extending
objects are enumerated by a polynomial in $k$. In order to find this extension it seems
easier to work with Gelfand-Tsetlin patterns rather than with strict plane partitions.
Next we define generalized Gelfand-Tsetlin patterns which 
turn out to be the right extension.

Let $r,n,c$ be integers, $r$ non-negative and  $n$ positive. In
this paper a {\it generalized $(r,n,c)$ Gelfand-Tsetlin-pattern}
(for short: $(r,n,c)$-pattern) is an array $(a_{i,j})_{1 \le i \le
r+1, i-1 \le j \le n+1}$ of integers with

\begin{enumerate}

\item $a_{i,i-1}=0$ and $a_{i,n+1}=c$,

\item if $a_{i,j} \le a_{i,j+1}$ then $a_{i,j} \le a_{i-1,j} \le a_{i,j+1}$

\item if $a_{i,j} > a_{i,j+1}$ then $a_{i,j} > a_{i-1,j} > a_{i,j+1}.$

\end{enumerate}

A $(3,6,c)$-pattern for example is of the form

\vspace{5mm}

\begin{tabular}{ccccccccccccccc}
 &   & & $0$ & & $a_{4,4}$ & & $a_{4,5}$ & & $a_{4,6}$ & & $c$ &  &  & \\
 & &$0$ &   & $a_{3,3}$ & & $a_{3,4}$ & & $a_{3,5}$ & & $a_{3,6}$ & & $c$ & & \\
 & $0$ & & $a_{2,2}$ & & $a_{2,3}$ & & $a_{2,4}$ & & $a_{2,5}$ & & $a_{2,6}$ & & $c$ & \\
$0$ & & $a_{1,1}$ & & $a_{1,2}$ & & $a_{1,3}$ & & $a_{1,4}$ & &
$a_{1,5}$ & & $a_{1,6}$ & & $c$,
\end{tabular}

\vspace{5mm}

\noindent such that every entry not in the top row is
between its northwest neighbour $w$ and its northeast neighbour
$e$, if $w \le e$ then weakly between, otherwise strictly between.
Thus

\vspace{5mm}

\begin{center}
\begin{tabular}{ccccccccccccccc}
 &   & & $0$ & & $3$ & & $-5$ & & $10$ & & $4$ &  &  & \\
 & &$0$ &   & $2$ & & $-2$ & & $3$ & & $8$ & & $4$ & & \\
 & $0$ & & $2$ & & $-1$ & & $2$ & & $4$ & & $7$ & & $4$ & \\
$0$ & & $0$ & & $0$ & & $1$ & & $2$ & & $5$ & & $6$ & & $4$
\end{tabular}
\end{center}

\vspace{5mm}

\noindent is an example of an $(3,6,4)$-pattern. Note that a
generalized $(n-1,n,c)$ Gelfand-Tsetlin-pattern $(a_{i,j})$ with
$0 \le a_{n,n} \le c$ is a Gelfand-Tsetlin-pattern with $n$ rows and parts in 
$\{0,1,\ldots,c\}$ as defined at the beginning of this section. This is because $0 \le a_{n,n} \le c$ implies that 
the third condition in the definition of a generalized Gelfand-Tsetlin-pattern 
never applies.

Next we introduce the sign of an $(r,n,c)$-pattern, since we actually have to work 
with a signed enumeration if $a_{n,n} \notin \{0,1,\ldots,c\}$.
A pair  $(a_{i,j}, a_{i,j+1})$ with
$a_{i,j} > a_{i,j+1}$ and $i \not= 1$ is called an {\it inversion}
of the $(r,n,c)$-pattern and $(-1)^{\# \, \text{of inversions}}$ is said to
be the {\it sign} of the pattern, denoted by $\sgn(a)$. The $(3,6,4)$-pattern in the example
above has altogether $6$ inversions and thus its sign is $1$.
We define the following
expression
$$
F(r,n,c;k_{1},k_{2},\dots,k_{n-r}) = \sum_{a} \sgn(a),
$$
where the sum is over all $(r,n,c)$-patterns $(a_{i,j})$ with
top row defined by $k_{i} = a_{r+1,r+i}$ for $i=1,\dots,n-r$. Now it is
important to observe that for $0 \le k \le c$ the number of
$(n-1,n,c)$-patterns with $a_{n,n}=k$ is given by $F(n-1,n,c;k)$. This is because an 
$(n-1,n,c)$-pattern with $0 \le a_{n,n} \le c$ has no inversions.
Thus $F(n-1,n,c;k)$ is the quantity we want to compute. It has the advantage
that it is defined for all integers $k$, whereas our original 
enumeration problem was only defined for $0 \le k \le c$.

\section{$F(n-1,n,c;k)$ is a polynomial in $k$ of degree $2n-2$ at most}
\label{polynomial}

In this section we establish Step~2 of the method above for our refinement of
the Bender-Knuth (ex-)Conjecture. The following recursion for $F(r,n,c;k_1,k_2,\ldots,k_{n-r})$ is
fundamental.
\begin{multline}
\label{rec} F(r,n,c;k_1,k_2,\ldots,k_{n-r})= \\
\sum_{l_1=0}^{k_1}
\sum_{l_2=k_1}^{k_2} \sum_{l_3=k_2}^{k_3} \ldots
\sum_{l_{n-r}=k_{n-r-1}}^{k_{n-r}} \sum_{l_{n-r+1}=k_{n-r}}^c
F(r-1,n,c;l_1,l_2,\ldots,l_{n-r+1}).
\end{multline}
It is obvious for $(k_{1},k_{2},\ldots,k_{n-r})$ with
$0 \le k_{1} \le k_{2} \le \ldots \le k_{n-r} \le c$. 
After recalling the extended definition of the summation symbol 
\eqref{sumext} one observes that the generalized $(r,n,c)$ Gelfand-Tstelin-patterns 
and $F(r,n,c;k_{1},\ldots,k_{n-r})$ were simply defined in such a way
that this recursion holds for arbitrary integer tuples
$(k_{1},\ldots,k_{n-r})$. 
This recursion together with the initial condition
$$
F(0,n,c;k_1,k_2,\ldots,k_n)=1
$$
implies the following lemma.

\begin{lem} Let $r,n$ be integers, $r$ non-negative and $n$
positive. Then $F(r,n,c;k_1,\ldots,k_{n-r})$ can be expressed by a
polynomial in the $k_i$'s and in $c$.
\end{lem}

In the following  $F(r,n,c;k_1,\ldots,k_{n-r})$ is identified with this
polynomial. In particular $F(n-1,n,c;k)$ is a polynomial $k$
and with this we have established the first half of Step~2 in our
method. Next we aim to show that
$F(r,n,c;k_1,\ldots,k_{n-r})$ is of degree $2r$ at most in every
$k_i$. This will imply that $F(n-1,n,c;k)$ is of degree $2n-2$ at
most  in $k$ and completes Step~2. However, this degree estimation
is complicated and takes Lemma~\ref{1} -- \ref{degree}.

The degree of $F(r,n,c;k_1,\ldots,k_{n-r})$ in $k_{i}$ is the degree
of 
\begin{equation}
\label{sum-kernel}
\sum_{l_i=k_{i-1}}^{k_i} \sum_{l_{i+1}=k_i}^{k_{i+1}}
F(r-1,n,c;l_1,\ldots,l_{n-r+1}),
\end{equation}
in $k_i$, where $k_0=0$ and $k_{n-r+1}=c$. If we assume  
by induction with respect to $r$ that the degree of 
$F(r-1,n,c;l_{1},\ldots,l_{n-r})$ in $l_{i}$ is at most $2r-2$ as well as 
the degree in $l_{i+1}$, this observation only allows us to conclude easily that the 
degree of $F(r,n,c;k_{1},\ldots,k_{n-r})$ in $k_{i}$ is at most
$4r-2$, however, we want to establish that the degree is at most $2r$.
The following lemma shows how to obtain a sharper 
degree estimation in summations of our type.

In order to state this lemma we have to define an operator $D_i$ which turns out to be crucial for
the analysis of the recursion in \eqref{rec}. Let
$G(k_1,k_2,\ldots,k_m)$ be a function in $m$ variables and $1 \le
i \le m-1$. We set
\begin{multline*}
D_i G(k_1,\ldots,k_m) := \\
G(k_1,\ldots,k_{i-1},k_i,k_{i+1},k_{i+2},\ldots,k_m) +
G(k_1,\ldots,k_{i-1},k_{i+1}+1,k_{i}-1,k_{i+2},\ldots,k_m).
\end{multline*}

\begin{lem}
\label{1} Let $F(x_1,x_2)$ be a polynomial in $x_1$ and $x_2$ which is in
$x_1$ as well as in $x_2$ of degree at most $R$. Moreover assume that
$D_{1} F(x_1,x_2)$ is of degree $R$ as a polynomial in $x_1$ and
$x_2$, i.e. a linear combination of monomials $x_1^m x_2^n$ with $m+n
\le R$. Then $\sum\limits_{x_1=a}^y \sum\limits_{x_2=y}^b F(x_1,x_2)$ is
of degree at most $R+2$ in $y$.
\end{lem}

{\it Proof.} Set $F_1(x_1,x_2)=D_{1} F(x_1,x_2)/2$
and $F_2(x_1,x_2)=(F(x_1,x_2)-F(x_2+1,x_1-1))/2$. Clearly
$F(x_1,x_2)=F_1(x_1,x_2)+F_2(x_1,x_2)$. Observe that $F_2(x_2+1,x_1-1)=-F_2(x_1,x_2)$.
Thus $F_2(x_1,x_2)$ is a linear combination of terms of the form $(x_1)_m (x_2+1)_n - (x_1)_n
(x_2+1)_m$ with $m,n \le R$. Now observe that
\begin{multline*}
\sum_{x_1=a}^y \sum_{x_2=y}^b (x_1)_m (x_2+1)_n - (x_1)_n (x_2+1)_m = \\
\frac{1}{m+1} \frac{1}{n+1}
\left( (a-1)_{n+1} (b+1)_{m+1} - (a-1)_{m+1} (b+1)_{n+1} - (a-1)_{n+1} (y)_{m+1} + \right. \\
\left. (b+1)_{n+1} (y)_{m+1} + (a-1)_{m+1} (y)_{n+1} - (b +
1)_{m+1} (y)_{n+1} \right).
\end{multline*}
and thus $\sum\limits_{x_1=a}^y \sum\limits_{x_2=y}^b F_2(x_1,x_2)$ is a
polynomial of degree at most $R+1$ in $y$. By the assumption in
the lemma $\sum\limits_{x_1=a}^y \sum\limits_{x_2=y}^b F_1(x_1,x_2)$ is of
degree at most $R+2$ in $y$ and the assertion follows. \qed

\medskip

Thus it suffices to show that $D_{i}
F(n,r,c;.)(k_{1},\ldots,k_{n-r})$ is of degree at most $2r$ as a polynomial in
$k_{i}$ and $k_{i+1}$. In Lemma~\ref{decomp} we show a much stronger 
assertion, namely we prove a formula which expresses 
$D_{i} F(r,n,c;k_{1},\ldots,k_{n-r})$ as a product 
of $F(r,n-2,c+2;k_{1},\ldots,k_{i-1},k_{i+2}+2,\ldots,k_{n-r}+2)$ and 
an (explicit) polynomial in $k_{i}$ and $k_{i+1}$ which is obviously 
of degree $2r$. For the proof of Lemma~\ref{decomp} we need two other 
Lemmas and the following definition.

The operator $\Phi_{m}$, applicable to functions in $m$
variables and related to the recursion in \eqref{rec}, is defined as follows.
$$
\Phi_{m} G (k_1,\ldots,k_{m+1}) = \sum_{l_1=k_1}^{k_2}
\sum_{l_2=k_2}^{k_3} \ldots \sum_{l_{m}=k_m}^{k_{m+1}}
G(l_1,\ldots,l_{m}).
$$
Observe that \eqref{rec} is equivalent to the following.
$$
F(r,n,c;k_1,\ldots,k_{n-r}) = \Phi_{n-r+1} F(r-1,n,c;.)
(0,k_1,\ldots,k_{n-r},c).
$$

\begin{lem}
\label{fund}
Let $m$ be a positive integer, $1 \le i \le m$ and
$G({\bf l})$ be a function in ${\bf l}=(l_{1},\ldots,l_{m})$. Then
\begin{multline*}
D_i \Phi_m G(k_{1},\ldots,k_{m+1}) \\ 
= - \frac{1}{2} \left(
\sum_{l_1=k_{1}}^{k_2} \ldots \sum_{l_{i-2}=k_{i-2}}^{k_{i-1}}
\sum_{l_{i-1}=k_i+1}^{k_{i+1}+1} \sum_{l_{i}=k_i}^{k_{i+1}}
\sum_{l_{i+1}=k_{i}-1}^{k_{i+2}} \sum_{l_{i+2}=k_{i+2}}^{k_{i+3}}
\ldots
\sum_{l_{m}=k_{m}}^{k_{m+1}} D_{i-1} G ({\bf l}) \right. \\
\left. + \sum_{l_1=k_{1}}^{k_2} \ldots
\sum_{l_{i-2}=k_{i-2}}^{k_{i-1}} \sum_{l_{i-1}=k_{i-1}}^{k_{i}} \sum_{l_{i}=k_{i}}^{k_{i+1}}
\sum_{l_{i+1}=k_{i}-1}^{k_{i+1}-1}
\sum_{l_{i+2}=k_{i+2}}^{k_{i+3}} \ldots
\sum_{l_{m}=k_{m}}^{k_{m+1}} D_{i} G ({\bf l}) \right).
\end{multline*}
(Set $D_{0} G({\bf l})=0$ and $D_{m} G({\bf l})=0$.)
\end{lem}

{\it Proof.} We set
$$
g(l_{i-1},l_{i},l_{i+1}) = 
\sum_{l_{1}=k_{1}}^{k_{2}} \ldots \sum_{l_{i-2}=k_{i-2}}^{k_{i-1}}
\sum_{l_{i+2}=k_{i+2}}^{k_{i+3}} \ldots \sum_{l_{m}=k_{m}}^{k_{m+1}} G(l_{1},\ldots,l_{m}).
$$
It suffices to show the following.
\begin{multline}
\label{toshow}
\sum_{l_{i-1}=k_{i-1}}^{k_{i}} \sum_{l_{i}=k_{i}}^{k_{i+1}}
\sum_{l_{i+1}=k_{i+1}}^{k_{i+2}} g(l_{i-1},l_{i},l_{i+1}) + 
\sum_{l_{i-1}=k_{i-1}}^{k_{i+1}+1} \sum_{l_{i}=k_{i+1}+1}^{k_{i}-1}
\sum_{l_{i+1}=k_{i}-1}^{k_{i+2}} g(l_{i-1},l_{i},l_{i+1}) \\
= - \frac{1}{2} \left( 
\sum_{l_{i-1}=k_{i}+1}^{k_{i+1}+1} \sum_{l_{i}=k_{i}}^{k_{i+1}} \sum_{l_{i+1}=k_{i}-1}^{k_{i+2}}
g(l_{i-1},l_{i},l_{i+1}) + g(l_{i}+1,l_{i-1}-1,l_{i+1}) \right. \\
+ \left. 
\sum_{l_{i-1}=k_{i-1}}^{k_{i}} \sum_{l_{i}=k_{i}}^{k_{i+1}}
\sum_{l_{i+1}=k_{i}-1}^{k_{i+1}-1} 
g(l_{i-1},l_{i},l_{i+1}) + g(l_{i-1},l_{i+1}+1,l_{i}-1) \right)
\end{multline}
By \eqref{sumext} the left-hand-side of this equation is equal to 
\begin{multline*}
\sum_{l_{i-1}=k_{i-1}}^{k_{i}} \sum_{l_{i}=k_{i}}^{k_{i+1}}
\sum_{l_{i+1}=k_{i+1}}^{k_{i+2}} g(l_{i-1},l_{i},l_{i+1}) - 
\sum_{l_{i-1}=k_{i-1}}^{k_{i+1}+1} \sum_{l_{i}=k_{i}}^{k_{i+1}}
\sum_{l_{i+1}=k_{i}-1}^{k_{i+2}} g(l_{i-1},l_{i},l_{i+1}) \\
=
\sum_{l_{i-1}=k_{i-1}}^{k_{i}} \sum_{l_{i}=k_{i}}^{k_{i+1}}
\sum_{l_{i+1}=k_{i+1}}^{k_{i+2}} g(l_{i-1},l_{i},l_{i+1}) - 
\sum_{l_{i-1}=k_{i-1}}^{k_{i}} \sum_{l_{i}=k_{i}}^{k_{i+1}}
\sum_{l_{i+1}=k_{i}-1}^{k_{i+2}} g(l_{i-1},l_{i},l_{i+1}) \\
- \sum_{l_{i-1}=k_{i}+1}^{k_{i+1}+1} \sum_{l_{i}=k_{i}}^{k_{i+1}}
\sum_{l_{i+1}=k_{i}-1}^{k_{i+2}} g(l_{i-1},l_{i},l_{i+1}) \\
= 
- \sum_{l_{i-1}=k_{i-1}}^{k_{i}} \sum_{l_{i}=k_{i}}^{k_{i+1}}
\sum_{l_{i+1}=k_{i}-1}^{k_{i+1}-1} g(l_{i-1},l_{i},l_{i+1}) - 
\sum_{l_{i-1}=k_{i}+1}^{k_{i+1}+1} \sum_{l_{i}=k_{i}}^{k_{i+1}}
\sum_{l_{i+1}=k_{i}-1}^{k_{i+2}} g(l_{i-1},l_{i},l_{i+1}).
\end{multline*}
The last expression is obviously equal to the right-hand-side of
\eqref{toshow} and the assertion of the lemma is proved. \qed

\medskip

\begin{lem}
\label{2}
Let $d$ and $r \ge 2$ be integers. Then
$$
\sum_{x'=x+d}^{y+d} \sum_{y'=x-1+d}^{y-1+d} (y'-x'-r+3)_{2r-3} (y'-x'+1) =
\frac{1}{r (2r-1)} (y-x-r+2)_{2r-1} (y-x+1).
$$
\end{lem}

{\it Proof.} Observe that
\begin{multline*}
\sum_{x'=x+d}^{y+d} \sum_{y'=x-1+d}^{y-1+d} (y'-x'-r+3)_{2r-3} (y'-x'+1) \\ =
\frac{1}{2} \sum_{x'=x+d}^{y+d} \sum_{y'=x-1+d}^{y-1+d} (y'-x'-r+2)_{2r-2} +
(y'-x'-r+3)_{2r-2}.
\end{multline*}
Furthermore
\begin{multline*}
\sum_{y'=x-1+d}^{y-1+d} (y'-x'-r+2)_{2r-2} \\ = \frac{1}{2r-1} \left(
  (y-x'-r+1+d)_{2r-1} - (x-x'-r+d)_{2r-1} \right) \\ =
\frac{1}{2r-1} \left(
  (x'-x-r+2-d)_{2r-1} - (x'-y-r+1-d)_{2r-1} \right)
\end{multline*}
and
\begin{multline*}
\sum_{y'=x-1+d}^{y-1+d} (y'-x'-r+3)_{2r-2} =
 \frac{1}{2r-1} \left(
  (x'-x-r+1-d)_{2r-1} - (x'-y-r-d)_{2r-1} \right).
\end{multline*}
Therefore the left-hand-side in the statement of lemma is equal to 
\begin{multline*}
\frac{1}{2 \, (2r-1)} \left ( \sum_{x'=x+d}^{y+d} (x'-x-r+2-d)_{2r-1} -
  (x'-y-r+1-d)_{2r-1}  \right. \\ \Bigg. + (x'-x-r+1-d)_{2r-1} -
  (x'-y-r-d)_{2r-1} \Bigg) \\ =
\frac{1}{4 r \, (2r-1)} \left( (y-x-r+2)_{2r} - (-r+1)_{2r} - (-r+1)_{2r} +
  (x-y-r)_{2r} \right. \\ \left. + (y-x-r+1)_{2r} - (-r)_{2r} - (-r)_{2r} + (x-y-r-1)_{2r} \right)
  \\  =
\frac{1}{2 r \, (2r-1)} \left( (y-x-r+2)_{2r} + (y-x-r+1)_{2r}\right),
\end{multline*}
since $(-r+1)_{2r}=0$ and $(-r)_{2r}=0$, and the assertion follows. \qed

\begin{lem}
\label{decomp}
Let $r,n,i$ be integers, $r$ non-negative, $n$ positive and $1
\le i \le n-r-1$. Then
\begin{multline*}
D_i F(r,n,c;.)(k_1,\ldots,k_{n-r}) = (-1)^{r}\frac{2}{(2r)!} (k_{i+1}-k_i-r+2)_{2r-1} (k_{i+1}-k_{i}+1) \\ \times F(r,n-2,c+2;k_{1},\ldots,k_{i-1},k_{i+2}+2,\ldots,k_{n-r}+2).
\end{multline*}
\end{lem}

{\it Proof.} We show the assertion by induction with respect to $r$. For $r=0$
there is nothing to prove. We assume $r > 0$. By the induction hypothesis we
may assume that
\begin{multline}
\label{i}
D_{i} F(r-1,n,c,.)(l_{1},\ldots,l_{n-r+1}) = (-1)^{r-1} \frac{2}{(2r-2)!}
(l_{i+1}-l_{i}-r+3)_{2r-3} (l_{i+1}-l_{i}+1)
\\ \times F(r-1,n-2,c+2;l_{1},\ldots,l_{i-1},l_{i+2}+2,\ldots,l_{n-r+1}+2)
\end{multline}
and
\begin{multline}
\label{i+1}
D_{i+1} F(r-1,n,c,.)(l_{1},\ldots,l_{n-r+1}) = (-1)^{r-1} \frac{2}{(2r-2)!} (l_{i+2}-l_{i+1}-r+3)_{2r-3}
(l_{i+2}-l_{i+1}+1) \\ \times F(r-1,n-2,c+2;l_{1},\ldots,l_{i},l_{i+3}+2,\ldots,l_{n-r+1}+2).
\end{multline}
By \eqref{rec} we have
\begin{equation*}
D_i F(r,n,c;.)(k_1,\ldots,k_{n-r}) =
D_{i} \Phi_{n-r+1} F(r-1,n,c;.)(0,k_{1},\ldots,k_{n-r},c).
\end{equation*}
Lemma~\ref{fund} implies that this is equal to
\begin{multline*}
-\frac{1}{2} \left(
\sum_{l_{1}=0}^{k_{1}} \ldots \sum_{l_{i-1}=k_{i-2}}^{k_{i-1}}
\sum_{l_{i}=k_{i}+1}^{k_{i+1}+1} \sum_{l_{i+1}=k_{i}}^{k_{i+1}}
\sum_{l_{i+2}=k_{i}-1}^{k_{i+2}} \ldots \sum_{l_{n-r+1}=k_{n-r}}^{c}
D_{i} F(r-1,n,c;.)(l_{1},\ldots,l_{n-r+1}) \right. \\
\left.+
\sum_{l_{1}=0}^{k_{1}} \ldots \sum_{l_{i}=k_{i-1}}^{k_{i}}
\sum_{l_{i+1}=k_{i}}^{k_{i+1}} \sum_{l_{i+2}=k_{i}-1}^{k_{i+1}-1}
\sum_{l_{i+3}=k_{i+2}}^{k_{i+3}} \ldots \sum_{l_{n-r+1}=k_{n-r}}^{c}
D_{i+1} F(r-1,n,c;.)(l_{1},\ldots,l_{n-r+1}) \right).
\end{multline*}
In this expression we replace $D_{i}
F(r-1,n,c;.)(l_{1},\ldots,l_{n-r+1})$ by the right-hand-side of
\eqref{i} and $D_{i+1} F(r-1,n,c;.)(l_{1},\ldots,l_{n-r+1})$ by
the right-hand-side of \eqref{i+1}. We note that by Lemma~\ref{2}
\begin{multline*}
\sum_{l_{i}=k_{i}+1}^{k_{i+1}+1} \sum_{l_{i+1}=k_{i}}^{k_{i+1}}
(l_{i+1}-l_{i}-r+3)_{2r-3} (l_{i+1}-l_{i}+1) \\ =
\sum_{l_{i+1}=k_{i}}^{k_{i+1}} \sum_{l_{i+2}=k_{i}-1}^{k_{i+1}-1}
(l_{i+2}-l_{i+1}-r+3)_{2r-3} (l_{i+2}-l_{i+1}+1) \\
= \frac{1}{r \, (2r-1)} (k_{i+1}-k_{i} - r +2)_{2r-1} (k_{i+1}-k_{i} +1).
\end{multline*}
Consequently we obtain the following for the left-hand-side in the statement
of the lemma.
\begin{multline*}
(-1)^{r}\frac{2}{(2r)!} (k_{i+1}-k_i-r+2)_{2r-1} (k_{i+1}-k_{i}+1) \\
\times \left( \sum_{l_{1}=0}^{k_{1}} \ldots \sum_{l_{i-1}=k_{i-2}}^{k_{i-1}}
\sum_{l_{i+2}=k_{i}-1}^{k_{i+2}} \sum_{l_{i+3}=k_{i+2}}^{k_{i+3}} \ldots \sum_{l_{n-r+1}=k_{n-r}}^{c} \right. \\
 F(r-1,n-2,c+2;l_{1},\ldots,l_{i-1},l_{i+2}+2,l_{i+3}+2,\ldots,l_{n-r+1}+2)   \\
+  \sum_{l_{1}=0}^{k_{1}} \ldots \sum_{l_{i-1}=k_{i-2}}^{k_{i-1}} \sum_{l_{i}=k_{i-1}}^{k_{i}}
\sum_{l_{i+3}=k_{i+2}}^{k_{i+3}} \ldots \sum_{l_{n-r+1}=k_{n-r}}^{c} \\
\Bigg.  F(r-1,n-2,c+2;l_{1},\ldots,l_{i-1},l_{i},l_{i+3}+2,\ldots,l_{n-r+1}+2)
\Bigg).
\end{multline*}
This is equal to
\begin{multline*}
(-1)^{r}\frac{2}{(2r)!} (k_{i+1}-k_i-r+2)_{2r-1} (k_{i+1}-k_{i}+1) \\
\times \sum_{l_{1}=0}^{k_{1}} \ldots
\sum_{l_{i}=k_{i-1}}^{k_{i+2}+2} \ldots
\sum_{l_{n-r+1}=k_{n-r}+2}^{c+2}  
F(r-1,n-2,c+2;l_{1},\ldots,l_{i-1},l_{i},l_{i+3},\ldots,l_{n-r+1}) \\
= (-1)^{r}\frac{2}{(2r)!} (k_{i+1}-k_i-r+2)_{2r-1} (k_{i+1}-k_{i}+1) \\
\times F(r,n-2,c+2;k_{1},\ldots,k_{i-1},k_{i+2}+2,\ldots,k_{n-r}+2)
\end{multline*}
and the assertion follows. \qed

\medskip

We are finally able to prove the degree lemma.

\begin{lem}
\label{degree} Let $r,n,i$ be integers, $r$ non-negative, $n$
positive and $1 \le i \le n-r$. Then $F(r,n,c;k_1,\ldots,k_{n-r})$
is a polynomial in $k_i$ of degree at most $2r$.
\end{lem}

{\it Proof.} We prove the assertion by induction with respect to
$r$. For $r=0$ it is trivial. Assume $r>0$ and $1 \le i \le n-r$.
The degree of $F(r,n,c;k_1,\ldots,k_{n-r})$ in
$k_i$ is at most the degree of \eqref{sum-kernel} in $k_{i}$.
By Lemma~\ref{decomp}
the degree of $D_i F(r-1,n,c;l_1,\ldots,l_{n-r+1})$ as a
polynomial in $l_i$ and $l_{i+1}$ is $2r-2$. Moreover the degree
of $F(r-1,n,c;l_1,\ldots,l_{n-r+1})$ in $l_i$ as well as in
$l_{i+1}$ is at most $2r-2$ by the induction hypothesis. The
assertion follows from Lemma~\ref{1}. \qed

\section{Exploring the zeros of $F(n-1,n,c;k)$}
\label{zeros}

We finally establish Step~3 of our method for the refinement of
the Bender-Knuth (ex-)Conjecture.

\begin{lem}
\label{factor}
Let $r,n,i$ be integers, $r$ non-negative, $n$
positive and $1 \le i \le n-r$. Then $F(r,n,c;k_1,\ldots,k_{n-r})$
vanishes for $k_1=-1,-2,\ldots,-r$ and
$k_{n-r}=c+1,c+2,\ldots,c+r$.
\end{lem}

{\it Proof.} It suffices to show that there exists no
$(r,n,c)$-pattern with first row 
$$(0,k_1,\ldots,k_{n-r},c),$$ 
if $k_1=-1,-2,\ldots,-r$ or $k_{n-r}=c+1,c+2,\ldots,c+r$. 
Suppose $(a_{i,j})$ is an
$(r,n,c)$-pattern with $a_{r+1,r+1} \in \{-1,-2,\ldots,-r\}$. In
particular we have $0 > a_{r+1,r+1}$ and thus the definition of
$(r,n,c)$-patterns implies that $0 > a_{r,r} > a_{r+1,r+1}$. In a
similar way we obtain $0 > a_{1,1} > a_{2,2} > \ldots > a_{r,r}
> a_{r+1,r+1}$. This is, however, a contradiction, since there
exist no $r$ distinct integers strictly between $0$ and $a_{r+1,r+1}$
if $a_{r+1,r+1} \in \{-1,-2,\ldots,-r\}$. The
case that $a_{n-r} \in \{c+1,c+2,\ldots,c+r\}$ is similar. \qed

\medskip

We obtain the following corollary.

\begin{cor}
\label{independent}
$F(n-1,n,c;k) / ((1+k)_{n-1} (1+c-k)_{n-1})$
is independent of $k$.
\end{cor}

{\it Proof.} By Lemma~\ref{factor} $(k+1)_{n-1} (1+c-k)_{n-1}$ is
a factor of $F(n-1,n,c;k)$. By Lemma~\ref{degree} $F(n-1,n,c;k)$
is of degree at most $2n-2$ in $k$ and the assertion follows. \qed

\begin{theo}
\label{special}
The number of strict plane partitions with parts
in $\{1,2,\ldots,n\}$, at most $c$ columns and $k$ parts equal to
$n$ is given by
$$
F(n-1,n,c;k) = \frac{(1+k)_{n-1} (1+c-k)_{n-1}}{(1)_{n-1}}
\prod_{i=1}^{n-1} \frac{(c+i+1)_{i-1}}{(i)_i}.
$$
\end{theo}

{\it Proof.} We prove the assertion by induction with respect to
$n$. Observe that the formula is true for $n=1$ since
$F(0,1,c;k)=1$. Assume $n > 1$. By Corollary~\ref{independent}
\begin{equation*}
F(n-1,n,c;k)= (1+k)_{n-1} (1+c-k)_{n-1}
\frac{F(n-1,n,c;c)}{(1+c)_{n-1} (1)_{n-1}}.
\end{equation*}
Observe that if we have $a_{n,n}=c$ in an $(n-1,n,c)$-pattern
$(a_{i,j})_{1 \le i \le n, i-1
\le j \le n+1}$ then $a_{i,n}=c$ for all $i$.
This implies the recursion
\begin{equation}
\label{extra} F(n-1,n,c;c)=\sum_{k=0}^c F(n-2,n-1,c;k).
\end{equation}
We need one other ingredient, namely the following hypergeometric
identity
\begin{multline}
\label{hyper} \sum_{k=0}^c (1+k)_{m-1} (1+c-k)_{m-1} 
= (1)^2_{m-1} \sum_{k=0}^c \binom{m+k-1}{m-1} \binom{c-k+m-1}{m-1}
\\
= (1)^2_{m-1} \binom{c+2m-1}{2m-1} 
= \frac{(1)^2_{m-1} (c+1)_{2m-2}}{(1)_{2m-2}},
\end{multline}
where the second equation is equivalent to the  Chu-Vandermonde identity; see
\cite[p. 169, (5.26)]{knuth}. With the help of the
recursion \eqref{extra}, the induction hypothesis for
$F(n-2,n-1,c;k)$ and the hypergeometric identity we are able to
compute $F(n-1,n,c;c)$ and with this $F(n-1,n,c;k)$. \qed

\begin{rem}
By the symmetry of the Schur function, the
number of strict plane partitions of a fixed shape with $x_i$
parts equal to $i$ is equal to the number of strict plane
partitions with $x_{\pi(i)}$ parts equal to $i$ for every
permutation $\pi$. Thus Theorem~\ref{special} gives the number of
strict plane partitions with parts in $\{1,2,\dots,n\}$, at most
$c$ columns and $k$ parts are equal to $i$ for arbitrary $i
\in \{1,2,\dots,n\}$. Also note that this does not generalize to
the generating function of these objects.
\end{rem}

\medskip

\begin{cor}[Andrews~\cite{andrews}, Gordon~\cite{gordon},
  Macdonald~\cite{macdonald}, Proctor~\cite{proctor}]
The number of strict plane partitions with parts in
$\{1,2,\dots,n\}$ and at most $c$ columns is
$$
\prod_{i=1}^n \frac{(c+i)_i}{(i)_i}.
$$
\end{cor}

{\it Proof.} By Theorem~\ref{special} the number of strict plane
partitions with parts in $\{1,2,\dots,n\}$ and at most $c$ columns
equals
$$
\sum_{k=0}^c \frac{(k+1)_{n-1} (1+c-k)_{n-1}}{(1)_{n-1}}
\prod_{i=1}^{n-1} \frac{(c+i+1)_{i-1}}{(i)_{i}}.
$$
The assertion now follows from (\ref{hyper}).

\section{Semistandard tableaux of a fixed shape}
\label{semistand}

In this section we apply our method to the enumeration of
semistandard tableaux of a fixed shape. This result is
certainly well-known. Nonetheless we think it might be interesting for the reader
to see another application of our method which moreover uses more
than just one ``polynomial parameter'' as opposed to the single parameter $k$
in the example above. (At this point the reader may wonder what we mean by a
multivariate application of our method, since we only describe the case of a single
polynomial parameter in the introduction. However, it is
straightforward to generalize this method to a multivariate version as it should become
clear in this section.)

Let $\lambda=(\lambda_{1},\lambda_{2},\dots,\lambda_{r})$ be a
partition and $k$ a positive integer. A semistandard tableau of
shape $\lambda$ with entries between $1$ and $k$ is a filling of
the Ferrers diagram of shape $\lambda$ with entries weakly between
$1$ and $k$ such that the rows are weakly increasing and the
columns are strictly increasing. (Semistandard tableaux and strict plane
partitions are equivalent objects. Indeed, if we replace every entry $e$ in a
semistandard tableau with entries between $1$ and $k$ with $1+k-e$ we clearly
obtain a strict plane partition. However, we choose to use the notion of semistandard
tableaux in this section for historical reasons.) It is well-known \cite[p. 375, in
(7.105) $q \to 1$]{stan} that the number of semistandard tableaux of
shape $\lambda$ with entries between $1$ and $k$ is
$$
\prod_{1 \le i < j \le r} \frac{\lambda_{i} - \lambda_{j} + j - i}{j-i}
\prod_{i=1}^{r} \frac{(\lambda_{i}+r+1-i)_{k-r}}{(r+1-i)_{k-r}}
$$
if $r \le k$, otherwise this number is obviously zero by the
columnstrictness. If $r=k$ the formula simplifies to
\begin{equation}
\label{formula} \prod_{1 \le i < j \le k} \frac{\lambda_{i} -
\lambda_{j} + j -i}{j-i}.
\end{equation}
It suffices to prove this formula, for if we have $r < k$ then 
the number of semistandard tableaux of shape 
$(\lambda_{1}, \ldots, \lambda_{r})$ is obviously equal to the number 
of semistandard tableaux of shape
$(\lambda_{1},\ldots,\lambda_{r},0,\ldots,0)$ ($k-r$ zeros).

The expression in (\ref{formula}) is a polynomial in the
$\lambda_{i}$'s which is up to a constant determined by its zeros
$\lambda_{i} =  \lambda_{j} - j + i$, $1 \le i < j \le k$. Clearly
the number of semistandard tableaux of shape
$\lambda=(\lambda_{1},\lambda_{2},\dots,\lambda_{k})$ with entries
between $1$ and $k$ can be interpreted to be zero if $\lambda_{i}
=  \lambda_{j} - j + i$ for some $i,j$ with $1 \le i < j \le k$,
since $\lambda$ is not a partition in this case. However, it is wrong to 
conclude that the appropriate combinatoral extension for the number of
semistandard tableaux of shape $(\lambda_{1},\ldots,\lambda_{k})$ is to set this number
to zero whenever $(\lambda_{1},\ldots,\lambda_{k})$ is not a
partition. In fact we will see that this number has to be zero if and
only if  $\lambda_{i}
=  \lambda_{j} - j + i$ for some $i,j$ with $1 \le i < j \le k$. Again we
divide the proof of \eqref{formula} into three steps.

\medskip

{\it Step 1.}
We extend the combinatorial interpretation of the number of semistandard tableaux of shape
$\lambda$ to arbitrary $\lambda \in \mathbb{Z}^{k}$. Define a function $F_k$
from $k$-tuples of integers to integers as follows.

\begin{enumerate}
\item[(i)] If $\lambda_{1} > \lambda_{2} > \ldots \lambda_{k} \ge -k$ then
$F_k(\lambda_{1},\ldots,\lambda_{k})$ is the number of semistandard tableaux
of shape $(\lambda_{1}+1,\lambda_{2}+2,\ldots,\lambda_{k}+k)$ with entries in $\{1,2,\ldots,k\}$.

\item[(ii)] $F_k$ is invariant under adding the same integer to all arguments. 

\item[(iii)] If one permutes the argument of $F_k$ by a permutation $\pi$, the effect 
is to multliply the value of $F_k$ by $\sgn(\pi)$, i.e.
\begin{equation}
\label{alternate}
F_k(\lambda_{1},\lambda_{2},\ldots,\lambda_{k}) = \sgn \pi \, 
F_k(\lambda_{\pi(1)},\lambda_{\pi(2)},\ldots,\lambda_{\pi(k)}).
\end{equation}
\end{enumerate}

Observe that (i) and (ii) are not contradictory. This is 
because the number of semistandard tableaux of shape
$(\lambda_{1},\ldots,\lambda_{k})$ with entries in $\{1,2,\ldots,k\}$
is equal to the number of semistandard
tableaux of shape $(\lambda_{1}+c,\lambda_{2}+c,\ldots,\lambda_{k}+c)$, $c \ge 0$, with entries in
$\{1,2,\ldots,k\}$, for in the latter case the first $c$ columns are equal to 
$(1,2,\ldots,k)^{t}$. Also note that (iii) implies  $F_k(\lambda_{1},\ldots,\lambda_{k})=0$ if
$\lambda_{i} = \lambda_{j}$ for distinct $i,j$. In order to prove 
\eqref{formula} we have to show that 
$$
F_k(\lambda_{1},\ldots,\lambda_{k}) = \prod_{1 \le i < j \le k}
\frac{\lambda_{i} - \lambda_{j}}{j-i}.
$$

\medskip

{\it Step 2.} Next we aim to show that the function $(\lambda_{1},\ldots,\lambda_{k-1})
\to F_k(\lambda_{1},\ldots,\lambda_{k-1},0)$  can be
expressed by a polynomial in $(\lambda_{1},\ldots,\lambda_{k-1})$ of degree at most 
$k-1$ in every $\lambda_{i}$ if $\lambda_{i} \ge 0$. (More general Step 3 will then imply
that the $(\lambda_{1},\ldots,\lambda_{k}) \to
F_{k}(\lambda_{1},\ldots,\lambda_{k})$ can be expressed by a polynomial if 
$(\lambda_{1},\ldots,\lambda_{k})$ is an arbitrary integer tupel.) This is again done by a recursion, however, 
in this case it requires some work to deduce it. 

If $\lambda=(\lambda_{1},\lambda_{2},\ldots,\lambda_{k})$ is a strict partition then the possible cells for the entry
$k$ in a semistandard tableau of shape $(\lambda_{1}+1,\lambda_{2}+2,\ldots,\lambda_{k}+k)$ with entries
between $1$ and $k$ are the cells $(i,j)$ with $j >
\lambda_{i+1}+i+1$. Moreover, by the columnstrictness, every cell in
the $k$-th row must contain the entry $k$. This implies the
following recursion
\begin{equation}
\label{recursion} F_k (\lambda_{1},\lambda_{2},\dots,\lambda_{k}) =
\sum_{\mu_{1}=\lambda_{2}+1}^{\lambda_{1}}
\sum_{\mu_{2}=\lambda_{3}+1}^{\lambda_{2}} \ldots
\sum_{\mu_{k-1}=\lambda_{k}+1}^{\lambda_{k-1}}
 F_k (\mu_{1},\mu_{2},\dots,\mu_{k-1})
\end{equation}
if $\lambda$ is a strict partition.

Let $\lambda \in \mathbb{Z}^{k}$ be with $\lambda_{1} \ge \lambda_{2}  \ge
\ldots \ge \lambda_{k}$. We show
\begin{equation}
\label{zero} \sum_{(\mu_{1},\ldots,\mu_{k-1}): \lambda_k+1 \le \mu_{i} \le
\lambda_{i}, \exists i': \mu_{i'} \le \lambda_{i'+1}} F_{k-1} (\mu_{1},\mu_{2},\dots,\mu_{k-1}) = 0,
\end{equation}
where the sum is over all $(\mu_{1},\ldots,\mu_{k-1})$, $\mu_{i}$ weakly
between $\lambda_{k}+1$ and $\lambda_{i}$ such that there
exists an $i'$, $1 \le i' \le k-1$, with $\mu_{i'} \le
\lambda_{i'+1}$.  For such a $(\mu_{1},\mu_{2},\dots,\mu_{k-1}) \in
\mathbb{Z}^{k-1}$ let $i'$ be minimal with the property that $\mu_{i'}
\le \lambda_{i'+1}$. Observe that $i' < k-1$ since $\lambda_{k}+1 \le
\mu_{k-1} \le \lambda_{k-1}$ by assumption. Then
$$
F_k (\mu_{1},\dots,\mu_{i'-1}, \mu_{i'},\mu_{i'+1}, \mu_{i'+2},
\dots,\mu_{k}) = - F_k (\mu_{1},\dots,\mu_{i'-1},
\mu_{i'+1},\mu_{i'}, \mu_{i'+2},\dots,\mu_{k})
$$
by (\ref{alternate}) if we set $\sigma = (i',i'+1)$. This induces
a sign-reversing involution on the set of summands since $\lambda_k+1  
\le \mu_{i'} \le \lambda_{i'+1}$, $\lambda_k+1 \le \mu_{i'+1} \le \lambda_{i'+1} \le
\lambda_{i'}$ and $\mu_{i'+1} \le \lambda_{i'+1}$. Now
(\ref{zero}) follows.

If we merge (\ref{recursion}) and (\ref{zero}), we obtain
\begin{equation}
\label{rec1} F_k(\lambda_1,\lambda_2,\dots,\lambda_k)=
\sum_{\mu_1=\lambda_k+1}^{\lambda_{1}}
\sum_{\mu_{2}=\lambda_{k}+1}^{\lambda_{2}} \ldots
\sum_{\mu_{k-1}=\lambda_k+1}^{\lambda_{k-1}} F_{k-1}
(\mu_1,\mu_2,\dots,\mu_{k-1})
\end{equation}
if $\lambda=(\lambda_{1},\lambda_{2},\ldots,\lambda_{k})$ is a strict partition. By
(ii) it is also true for strictly decreasing integers sequences $\lambda$. 
Moreover it is easily extendible to weakly decreasing $\lambda \in
\mathbb{Z}^{k}$: If there exists an $i'$ with $\lambda_{i'}
= \lambda_{i'+1}$ the left hand side vanishes by definition. The right-hand-side
vanishes, since it is equal the left-hand side of (\ref{zero}), because every 
$(\mu_{1},\ldots,\mu_{k-1})$ in the summation domain satisfies $\mu_{i'} \le
\lambda_{i'} = \lambda_{i'+1}$.

Finally we extend \eqref{rec1} to $\lambda \in \mathbb{Z}^k$ with $\lambda_{i} \ge
\lambda_{k}$. In this case  there exists a permutation $\pi \in {\mathcal S}_k$ with
$\lambda_{\pi(1)} \ge \lambda_{\pi(2)} \ge \ldots \ge \lambda_{\pi(k)}$. Clearly
$\pi(k)=k$. Consequently
\begin{multline}
\label{id}
F_k (\lambda_{1},\lambda_{2},\dots,\lambda_{k}) = \sgn(\pi) \, F_k(\lambda_{\pi(1)},\ldots,\lambda_{\pi(k)}) \\
= \sgn \pi 
\sum_{\mu_1=\lambda_k+1}^{\lambda_{\pi(1)}}
\sum_{\mu_{2}=\lambda_{k}+1}^{\lambda_{\pi(2)}} \ldots
\sum_{\mu_{k-1}=\lambda_k+1}^{\lambda_{\pi(k-1)}} F_{k-1}
(\mu_1,\mu_2,\dots,\mu_{k-1}) \\
= \sgn \pi \sum_{\mu_{\pi(1)}=\lambda_k+1}^{\lambda_{\pi(1)}}
\sum_{\mu_{\pi(2)}=\lambda_{k}+1}^{\lambda_{\pi(2)}} \ldots
\sum_{\mu_{\pi(k-1)}=\lambda_k+1}^{\lambda_{\pi(k-1)}} 
F_{k-1}(\mu_{\pi(1)},\mu_{\pi(2)},\dots,\mu_{\pi(k-1)})  \\
=\sum_{\mu_1=\lambda_k+1}^{\lambda_{1}}
\sum_{\mu_{2}=\lambda_{k}+1}^{\lambda_{2}} \ldots
\sum_{\mu_{k-1}=\lambda_k+1}^{\lambda_{k-1}} F_{k-1}
(\mu_1,\mu_2,\dots,\mu_{k-1})
\end{multline}
Define 
$$
P_k(\lambda_{1},\lambda_{2},\ldots,\lambda_{k-1}) = 
\sum_{\mu_1=1}^{\lambda_{1}}
\sum_{\mu_{2}=1}^{\lambda_{2}} \ldots
\sum_{\mu_{k-1}=1}^{\lambda_{k-1}} F_{k-1}
(\mu_1,\mu_2,\dots,\mu_{k-1}).
$$
By induction with respect to $k$ (note that $F_{1}(\lambda_{1})=1$) we may assume that 
$F_{k-1}(\mu_{1},\ldots,\mu_{k-1})$ is a polynomial in the $\mu_{i}$'s which is of
degree $k-2$ at most in $\mu_{i}$.
Thus $P_k(\lambda_{1},\ldots,\lambda_{k-1})$ is a
polynomial of degree at most $k-1$ in every $\lambda_{i}$. By \eqref{id} it is equal to 
$F_k(\lambda_{1},\ldots,\lambda_{k-1},0)$ for $\lambda_{i} \ge 0$. The assertion
of this step follows.

\medskip

{\it Step 3.}
We know that $F_k(\lambda_{1},\ldots,\lambda_{k-1},0)$ vanishes if 
$\lambda_{i} = \lambda_{j}$, $1 \le i < j \le k-1$, or $\lambda_{i}=0$. Consequently
$P_k(\lambda_{1},\ldots,\lambda_{k})$
has the factors
$$
\prod_{1 \le i < j \le k-1} (\lambda_i - \lambda_j) 
\prod_{i=1}^{k-1} \lambda_i.
$$
(Note that we use the fact that a polynomial $Q(x_{1},\ldots,x_{n})$ which
vanishes for all $(x_{1},\ldots,x_{n}) \in \mathbb{Z}^{n}$ with $x_{i} \ge 0$
can only be the zero polynomial.)
The product of these two factors is a polynomial in the
$\lambda_i$'s, $i=1,2,\dots,k-1$, which is of degree $k-1$ in
$\lambda_i$. This determines $P_k(\lambda_{1},\ldots,\lambda_{k-1})$ up to
a constant. Since $P_k(k-1,k-2,\ldots,1)=F_k(k-1,k-2,\ldots,1,0)=1$ we have
$$
P_k(\lambda_1,\lambda_2,\dots,\lambda_{k-1}) = \prod_{1 \le i < j
\le k-1} \frac{\lambda_i - \lambda_j}{j-i} \,
\prod_{i=1}^{k-1} \frac{\lambda_i}{k-i}.
$$
Observe that by (ii) 
\begin{multline*}
F_k(\lambda_1,\lambda_2,\dots,\lambda_k)=F_k(\lambda_1-\lambda_k,\lambda_2-\lambda_k,\dots,
\lambda_{k-1}-\lambda_k,0)= \\ P_k(\lambda_1 - \lambda_k,\lambda_2 - \lambda_k,\dots,\lambda_{k-1}-\lambda_k)
\end{multline*}
if $\lambda_{i} \ge \lambda_{k}$ for all $i$. Consequently
\begin{equation}
\label{part} F_k(\lambda_1,\lambda_2,\dots,\lambda_k)
 = \prod_{1 \le i < j \le k}
\frac{\lambda_i - \lambda_j}{j-i}
\end{equation}
if $\lambda_{i} \ge \lambda_{k}$.

\medskip

Finally we have to show that \eqref{part} is valid for all
$\lambda \in \mathbb{Z}^{k}$. Again let $\pi \in {\mathcal S}_{k}$
be such that $\lambda_{\pi(1)}  \ge \lambda_{\pi(2)} \ge \ldots \ge \lambda_{\pi(k)}$. Then
\begin{multline*}
F_k(\lambda_{1},\lambda_{2},\dots,\lambda_{k}) 
= \sgn \pi \, F_k(\lambda_{\pi(1)},\lambda_{\pi(2)},\dots,\lambda_{\pi(k)})   \\
= \sgn \pi  \prod_{1 \le i < j \le k}
\frac{\lambda_{\pi(i)} - \lambda_{\pi(j)}}{j - i}  
=  \prod_{1 \le i < j \le k} \frac{\lambda_i - \lambda_j}{j-i}.
\end{multline*}
This concludes the proof of \eqref{formula}.

\begin{rem}
The extension of our method introduced in the following section can be used to
derive the $q$-version of \eqref{formula}, see \cite[p. 375, (7.105)]{stan}.
\end{rem}

\section{Extension of the method to $q$-polynomials}
\label{q-p}

A natural question to ask is whether it is possible to obtain a
generating function version of Theorem~\ref{special}. Of course
only this would refine the Bender-Knuth (ex-)Conjecture. Clearly
our generating function (see Theorem~\ref{main}) is not a polynomial in $k$, however, we
introduce the notion of a {\it $q$-polynomial} below and find that
the generating function is such a $q$-polynomial. Thus we
adapt our method to $q$-polynomials in this section.

\medskip

Let $I$ be an integral domain containing $\mathbb{Q}$. A {\it $q$-polynomial} over $I$ in
the variables $X_{1},X_{2},\ldots,X_{n}$ is an ordinary polynomial over
$I(q)$, the field of rational functions in $q$ over $I$, in $q^{X_{1}},
q^{X_{2}}, \ldots, q^{X_{n}}$. The ring of these $q$-polynomials is denoted by 
$I_{q}[X_{1},\ldots,X_{n}]$.
For expressions of the form 
$$
(q^{X_{1}})^{a_{1}} (q^{X_{2}})^{a_{2}} \ldots (q^{X_{n}})^{a_{n}} q^{c}
$$
in a $q$-polynomial, where the $a_{i}$ are integers, we also write 
$$
q^{a_{0}+a_{1} X_{1} + a_{2} X_{2} + \ldots + a_{n} X_{n}}.
$$
We define $[X;q]=(1-q^{X})/(1-q)$ and 
$[X;q]_{n}=\prod_{i=0}^{n-1} [X+i;q]$. Observe that 
$$
[X_{1};q]_{m_{1}} [X_{2};q]_{m_{2}} \ldots [X_{n};q]_{m_{n}},
$$
$(m_{1},m_{2},\ldots,m_{n}) \in \mathbb{Z}, m_{i} \ge 0$, is a basis of
$I_{q}[X_{1},\ldots,X_{n}]$ over $I(q)$. This basis is the most convenient for our 
purpose.

If we review the proof of Theorem~\ref{special} we see that the
following two basic properties of polynomials were crucial.

\begin{itemize}

\item If $p(X)$ is a polynomial over an integral domain containing $\mathbb{Q}$, then there exists a (unique)
polynomial $r(X)$ with $\deg r = \deg p +1$ and
$$
\sum_{x=1}^y p(x) = r(y)
$$
for every integer $y$.

\item If $p(X)$ is a polynomial over an integral domain containing
  $\mathbb{Q}$ and $a_1,a_2,\dots,a_r$ are distinct
zeros of $p(X)$, then there exists a polynomial $r(X)$ with
$$
p(X) = (X-a_1) (X-a_2) \dots (X-a_r) r(X).
$$

\end{itemize}

\medskip

The following analogous hold for $q$-polynomials.

\begin{itemize}

\item If $p(X)$ is a $q$-polynomial, then there exists a (unique)
$q$-polynomial $r(X)$ with $\deg r = \deg p + 1$ and
$$
\sum_{x=1}^y p(x) \, q^x = r(y)
$$
for all integers $y$. (The degree of a $q$-polynomial in $X$ is defined as 
the degree of the corresponding ordinary polynomial in $q^{X}$.) In order to see that note
$$
[X;q]_{n+1} - [X-1;q]_{n+1} = q^{X-1} [n+1;q] [X;q]_n,
$$
which implies
\begin{equation}
\label{q-poch}
\sum_{x=1}^y [x;q]_n q^x = \frac{q}{[n+1;q]} [y;q]_{n+1}
\end{equation}
for all integers $y$.

\item If $p(X)$ is a $q$-polynomial and $a_1,a_2,\dots,a_r$ are
distinct integer zeros of $p(X)$, then there exists a $q$-polynomial
$r(X)$ with
\begin{multline*}
p(X) = ([X;q] - [a_1;q]) ([X;q] - [a_2;q]) \cdots ([X;q] - [a_r;q]) r(X) = \\
q^{a_1+a_2+\ldots+a_r} [X - a_1;q] [X - a_2;q] \cdots [X - a_r;q]
r(X).
\end{multline*}
The proof is analogous to the proof for ordinary polynomials, namely
the fundamental identity is
$$
[X;q]^n-[a;q]^n = ([X;q] - [a;q]) \sum_{i=0}^{n-1} [X;q]^i
[a;q]^{n-1-i} = q^a [X - a;q] \sum_{i=0}^{n-1} [X;q]^i [a;q]^{n-1-i}.
$$

\end{itemize}

Using these $q$-analogs it is quite
straightforward to modifiy the proof of Theorem~\ref{special} in order to
prove Thorem~\ref{main}. In the following we sketch it by stating the
$q$-versions of the definitions and lemmas that were necessary to prove
Theorem~\ref{special}.

The norm of an $(r,n,c)$-pattern is defined as the sum of its parts, where we
omit the first and the last part in each row. Our first observation is that
the bijection in Lemma~\ref{bijection} is norm-preserving. We introduce a
$q$-analog of $F(r,n,c;k_{1},\ldots,k_{n-r})$. Let 
$$
F_{q}(r,n,c;k_{1},\ldots,k_{n-r}) = \left( \sum_{a} \sgn (a) q^{\norm (a)}
\right) / (q^{k_{1}+k_{2}+\ldots+k_{n-r}}),
$$
where the sum is over all $(r,n,c)$-patterns $(a_{i,j})$ with
$k_{i}=a_{r+1,r+i}$ for $i=1,2,\ldots,n-r$. Observe that 
$F(n-1,n,c;k) \, q^{k}$ is the generating function of strict plane 
partitions with parts in $\{1,2,\ldots,n\}$, at most $c$ columns 
and $k$ parts equal to $n$. We have
$F_{q}(0,n,c;k_{1},\ldots,k_{n})=1$
and 
\begin{multline}
\label{rec-q} F_{q}(r,n,c;k_1,k_2,\ldots,k_{n-r})= \\
\sum_{l_1=0}^{k_1}
\sum_{l_2=k_1}^{k_2} \sum_{l_3=k_2}^{k_3} \ldots
\sum_{l_{n-r}=k_{n-r-1}}^{k_{n-r}} \sum_{l_{n-r+1}=k_{n-r}}^c
F_{q}(r-1,n,c;l_1,l_2,\ldots,l_{n-r+1}) \, q^{l_{1}+l_{2}+\ldots+l_{n-r+1}}.
\end{multline}
This shows that $F_{q}(r,n,c;k_{1},k_{2},\ldots,k_{n-r})$ is a $q$-polynomial
in $(k_{1},\ldots,k_{n-r})$. Next we have to show that
$F_{q}(r,n,c;k_{1},\ldots,k_{n-r})$ is of degree $2r$ in $k_{i}$ at most.
For that propose we need the following $q$-analog of Lemma~\ref{1}.

\begin{lem}
\label{1-q} Let $F(x_{1},x_{2})$ be a $q$-polynomial in $x_{1}$ and $x_{2}$ which is in
$x_{1}$ as well as in $x_{2}$ of degree $R$ at most. Moreover assume that
$D_{1} F(x_{1},x_{2})$ is of degree $R$ at most as a $q$-polynomial in $x_{1}$ and
$x_{2}$, i.e. a linear combination of monomials $(q^{x_{1}})^m (q^{x_{2}})^n$ with $m+n
\le R$. Then $\sum\limits_{x_{1}=a}^y \sum\limits_{x_{2}=y}^b F(x_{1},x_{2}) q^{x_{1}+x_{2}}$ is
of degree $R+2$ at most in $y$.
\end{lem}

The $q$-version of the operator $\Phi_{m}$ is
defined as follows.
$$
\Phi^{q}_{m} G (k_1,\ldots,k_{m+1}) = \sum_{l_1=k_1}^{k_2}
\sum_{l_2=k_2}^{k_3} \ldots \sum_{l_{m}=k_m}^{k_{m+1}}
G(l_1,\ldots,l_{m}) q^{l_{1}+l_{2}+\ldots+l_{m}}.
$$
With this definition we are able to state the $q$-analog of Lemma~\ref{fund}.
\begin{lem}
\label{fund-q}
Let $m$ be a positive integer, $1 \le i \le m$ and
$G({\bf l})$ be a function in ${\bf l}=(l_{1},\ldots,l_{m})$. Then
\begin{multline*}
D_i \Phi^{q}_m G(k_{1},k_{2},\ldots,k_{m+1}) \\
= - \frac{1}{2} \left(
\sum_{l_1=k_{1}}^{k_2} \ldots \sum_{l_{i-2}=k_{i-2}}^{k_{i-1}}
\sum_{l_{i-1}=k_i+1}^{k_{i+1}+1} \sum_{l_{i}=k_i}^{k_{i+1}}
\sum_{l_{i+1}=k_{i}-1}^{k_{i+2}} \sum_{l_{i+2}=k_{i+2}}^{k_{i+3}}
\ldots
\sum_{l_{m}=k_{m}}^{k_{m+1}} D_{i-1} G ({\bf l}) \, q^{l_{1}+\ldots+l_{m}}\right. \\
\left. + \sum_{l_1=k_{1}}^{k_2} \ldots \sum_{l_{i-2}=k_{i-2}}^{k_{i-1}}
\sum_{l_{i-1}=k_{i-1}}^{k_{i}} \sum_{l_{i}=k_{i}}^{k_{i+1}}
\sum_{l_{i+1}=k_{i}-1}^{k_{i+1}-1}
\sum_{l_{i+2}=k_{i+2}}^{k_{i+3}} \ldots
\sum_{l_{m}=k_{m}}^{k_{m+1}} D_{i} G ({\bf l})  \, q^{l_{1}+\ldots+l_{m}} \right).
\end{multline*}
\end{lem}
It can be deduced from Lemma~\ref{fund} by applying it to
$G(l_{1},\ldots,l_{m}) \, q^{l_{1}+\ldots+l_{m}}$ rather than 
to $G(l_{1},\ldots,l_{m})$.
Next we state the $q$-analog of Lemma~\ref{2}.

\begin{lem}
\label{2-q}
 Let $d$ and $r \ge 2$ be integers. Then
\begin{multline*}
\sum_{x'=x+d}^{y+d} \sum_{y'=x-1+d}^{y-1+d} [y'-x'-r+3;q]_{2r-3}
  [y'-x'+1;q] \, q^{(2r-2)x'} \, (1+q^{r-1}) \, q^{x'+y'}  = 
\\
=
2  \frac{[y-x-r+2;q]_{2r-1}
  [y-x+1;q] q^{2r x}}{[2r-1;q] [2r;q]} (1+q^{r}) q^{2 \, d \, r + r -2}.
\end{multline*}
\end{lem}

{\it Proof.} The proof is analogous to the proof of Lemma~\ref{2}. The
fundamental identities are
$$
\sum_{z=a}^{b} [z+w]_{n} q^{z} = \frac{q^{-w+1}}{[n+1;q]} \left( [b+w;q]_{n+1}
  - [a-1+w;q]_{n+1} \right)
$$
which is an easy consequence of \eqref{q-poch},  and 
$$
[z;q]_{n} = (-1)^{n} q^{n(z+(n-1)/2)} [-z-n+1;q]_{n}. \qed
$$

Lemma~\ref{fund-q} and Lemma~\ref{2-q} imply the $q$-analog of Lemma~\ref{decomp}.

\begin{lem}
\label{decomp-q}
Let $r,n,i$ be integers, $r$ non-negative, $n$ positive and $1
\le i \le n-r-1$. Then
\begin{multline*}
D_i F_{q}(r,n,c;.)(k_1,\ldots,k_{n-r}) = (-1)^{r}\frac{(1+q^{r})}{[1;q]_{2r}}
[k_{i+1}-k_i-r+2;q]_{2r-1} [k_{i+1}-k_{i}+1;q] q^{2r k_{i}} \\ \times
q^{r(1+4i-4n+5r)/2} F_{q}(r,n-2,c+2;k_{1},\ldots,k_{i-1},k_{i+2}+2,\ldots,k_{n-r}+2).
\end{multline*}
\end{lem}

Lemma~\ref{decomp} shows
that $D_i F_{q}(r,n,c;.)(k_1,\ldots,k_{n-r})$ is of degree $2r$ in
$k_{i}$ and in $k_{i+1}$. In the next lemma we see that this is also true for $F(r,n,c;k_{1},\ldots,k_{n-r})$
itself.

\begin{lem}
\label{degree-q} Let $r,n,i$ be integers, $r$ non-negative, $n$
positive and $1 \le i \le n-r$. Then $F_{q}(r,n,c;k_1,\ldots,k_{n-r})$
is a $q$-polynomial in $k_i$ of degree at most $2r$.
\end{lem} 

{\it Proof.} Use \eqref{rec-q}, Lemma~\ref{decomp-q} and Lemma~\ref{1-q} in the
same way as their analogs in Lemma~\ref{degree}. \qed

\medskip

Next we state the $q$-analog of  Lemma~\ref{factor}, which deals with the zeros of 
$$F_{q}(r,n,c;k_{1},k_{2},\ldots,k_{n-r})$$ in $k_{1}$ and $k_{n-r}$.

\begin{lem}
\label{factor-q}
Let $r,n,i$ be integers, $r$ non-negative, $n$
positive and $1 \le i \le n-r$. Then $F_{q}(r,n,c;k_1,\ldots,k_{n-r})$
is zero for $k_1=-1,-2,\ldots,-r$ and for
$k_{n-r}=c+1,c+2,\ldots,c+r$.
\end{lem}

{\it Proof.} In the proof of Lemma~\ref{factor} we have showed that there
exists no $(r,n,c)$-pattern with first row $(0,k_{1},\ldots,k_{n-r},c)$ 
if $k_{1}=-1,-2,\ldots,-r$ or $k_{n-r}=c+1,c+2,\ldots,c+r$. \qed

\medskip

This, the previous lemma and the second property of $q$-polynomials imply the following $q$-analog of 
Corollary~\ref{independent}.

\begin{cor}
\label{q-ind}
$F_q(n-1,n,c;k) / ([1+k;q]_{n-1} [k-c-n+1;q]_{n-1})$ is
independent of $k$.
\end{cor}

Note that $[1+c-k;q]_{n-1}$ is not a $q$-polynomial in $k$ and
therefore we work with  $[k-c-n+1;q]_{n-1}$ instead. We are now able to prove our main theorem.

\medskip

{\it Proof of Theorem~\ref{main}.} We prove the assertion by induction with respect to
$n$. Observe that the formula is true for $n=1$ since
$F(0,1,c;k)=1$. Applying Corollary~\ref{q-ind} in the same way as 
Corollary~\ref{independent} was applied in the proof of Theorem~\ref{special},
it suffices to check the formula for $F_{q}(n-1,n,c;c)$.
For that purpose we need the recursion
$$
F_q(n-1,n,c;c) = q^{c \, n - c} \sum_{k=0}^{c} F_q(n-2,n-1,c;k) q^{k}
$$
and the following identity
\begin{equation}
\label{q-vand} \sum_{k=0}^c [k+1;q] _{m-1} [k-c-m+1;q]_{m-1} \,
q^k = (-1)^{m-1} q^{(-m+1)(2c+m)/2} \frac{[1;q]_{m-1}^2
[c+1;q]_{2m-1}}{[1;q]_{2m-1}}
\end{equation}
which can be deduced from the $q$-Chu-Vandermonde identity,
see~\cite[(1.5.3);Appendix (II. 6)]{gasper}. \qed

\medskip

Finally we are able to prove the Bender-Knuth (ex-)Conjecture.
\begin{cor}
The generating function of strict plane partitions with parts in
$\{1,2,\dots,n\}$ and at most $c$ columns is
$$
\prod_{i=1}^n \frac{[c+i;q]_i}{[i;q]_i}.
$$
\end{cor}

{\it Proof.} By Theorem~\ref{main} the generating function is
equal to
$$
\sum_{k=0}^c \frac{q^{k \, n} [k+1;q]_{n-1} [1+c-k;q]_{n-1}
}{[1;q]_{n-1}} \prod_{i=1}^{n-1}
\frac{[c+i+1;q]_{i-1}}{[i;q]_{i}}.
$$
The assertion follows from (\ref{q-vand}). \qed

\section{A final observation}
\label{final}

A monotone triangle of size $n$, see \cite[p. 58]{bressoud}, is an
$(n-1,n,n+1)$-pattern with strictly increasing rows.  Monotone triangles of size $n$ with the central part of
the first row equal to $k$ are easily seen to be in bijection with
alternating sign matrices of size $n$, where the unique $1$ in the
first row is in the $k$-th column. Let $A(n,k)$ denote the number
of these objects. It was conjectured by Mills, Robbins and Rumsey
\cite{mills} (well-known as the refined alternating sign matrix Theorem)
and proved by Zeilberger \cite{zeil} that
$$
A(n,k) = \frac{(k)_{n-1} (1+n-k)_{n-1}}{(1)_{n-1}} \prod_{i=1}^{n-1}
\frac{(1)_{3i-2}}{(1)_{n+i-1}}.
$$
Surprisingly it turns out that the number of
$(n-1,n,n-1)$-patterns $(a_{i,j})$ with $a_{n,n}=k-1$ divided by $A(n,k)$   is
independent of $k$. In fact it is equal to
$$
\prod_{1 \le i \le j \le n-1} \frac{i + j+ n - 2}{i+2 j - 2},
$$
the number of $(n-1) \times (n-1) \times (n-1)$ totally symmetric
plane partitions, see~\cite{stem}. In a similar way as for the enumeration
of $(n-1,n,c)$-patterns, it suffices to show that
$$A(n,k) / ((k)_{n-1} (1+n-k)_{n-1})$$
is independent of $k$ in order to prove the formula for $A(n,k)$,
see \cite[Sec. 5.2]{bressoud} for an explanation. Therefore we
hope to find another proof of the refined alternating sign matrix
Theorem which is along the lines of the proof of
Theorem~\ref{special}. The situation is similar to  the strict
plane partitions which are under consideration in this paper:
First, one has to find an extension of the combinatorial
interpretation of alternating sign matrices of order $n$ with the
unique $1$ in the first row is in the $k$-th column to {\it
arbitrary} integers $k$. That is to say that one has to find
combinatorial objects indexed by a positive integer $n$ and an
arbitrary integer $k$ which are in bijection with alternating sign
matrices of order $n$, where the unique $1$ in the first row is in
the $k$-th column for $k \in \{1,2,\dots,n\}$. In the view of fact
that the generalized $(n-1,n,c)$ Gelfand-Tsetlin-patterns were the
right extension of the strict plane partitions, one would rather
work with monotone triangles than with alternating sign matrices.
Next it has to be shown that for fixed $n$ these objects are
enumerated by a polynomial $P_n(k)$ in $k$ of degree $2n - 2$,
typically this could be done by a recursion similar to the one
given  in Lemma~\ref{2}. Finally it has to be shown that there
exist none of these extending combinatorial objects  if
$k=0,-1,\dots,-n+2$ or $k=n+1,n+2,\dots,2n-1$.

Finally observe the following: We have seen that in order to give
another proof of the refined alternating sign matrix
(ex-)Conjecture it would suffice to show that the number of
$(n-1,n,n-1)$-patterns $(a_{i,j})$ with $a_{n,n}=k-1$ divided by
the number of alternating sign matrices of order $n$, where the
unique $1$ in the first row is in the $k$-th column is independent
of $k$. Thus, a bijection between $(n-1,n,n-1)$-patterns with
$0 \le a_{n,n}=k-1 \le n-1$ on one side and pairs consisting of a
monotone triangle of size $n$ with the central part in the first
row equal to $k$ and $(n-1) \times (n-1) \times (n-1)$ totally
symmetric plane partitions would simultaneously prove the formula
for $A(n,k)$ and for the number of  $(n-1) \times (n-1) \times
(n-1)$ totally symmetric plane partitions.


\medskip







\begin{thebibliography} {6}

\bibitem{andrews}
G. E. Andrews, Plane Partitions II: The equivalence of the Bender-Knuth and the
MacMahon conjectures, {\it Pacific J. Math} {\bf 72} (1977), no. 2,  283 -- 291.

\bibitem{benderknuth}
E. A. Bender and D. E. Knuth, Enumeration of Plane Partitions,
{\it J. Combin. Theory Ser. A} { \bf 13} (1972), 40 -- 54.

\bibitem{bressoud}
D. M. Bressoud, {\it Proof and Confirmations, The Story of the Alternating Sign Matrix Conjecture},
Cambridge University Press, Cambridge, 1999.

\bibitem{des1}
J. D\'esarm\'enien, La d\'emonstration des identit\'es de Gordon et MacMahon et de deux identit\'es nouvelles, in: Actes de $15^\text{e}$ Seminaire Lotharingien, I.R.M.A. Strasbourg, 1987, 39 -- 49.

\bibitem{des2}
J. D\'esarm\'enien, Une g\'en\'eralisation des formules de Gordon et de MacMahon, {\it Comptes Rendus Acad. Sci. Paris, S\'erie I} {\bf 309}, (1989), no. 6, 269 -- 272.

\bibitem{gasper}
G. Gasper and R. Rahman, {\it Basic hypergeometric series}, Encyclopedia of Mathematics and its
Applications {\bf 35}, Cambridge University Press, Cambridge, 1990.

\bibitem{gel}
I. M. Gelfand and M. L. Tsetlin, Finite-dimensional representations of the group of unimodular
matrices (in Russian), {\it Doklady Akad. Nauk. SSSR (N. S.)} {\bf 71} (1950), 825 -- 828.

\bibitem{gordon}
B. Gordon, A proof of the Bender-Knuth Conjecture, {\it Pacific J.
Math.} {\bf 108} (1983), no. 1,  99 -- 113.

\bibitem{kadell}
K. W. J. Kadell, Sch\"utzenberger's jeu de taquin and plane partitions, {\it J. Combin. Theory Ser. A} {\bf 77} (1997), no. 1, 110 -- 133.

\bibitem{knuth}
R. L. Graham, D. E. Knuth, O. Patashnik, {\it Concrete Mathematics }, 
2nd ed., Addison-Wesley, Reading, MA, 1994.

\bibitem{kratt}
C. Krattenthaler, The major counting of nonintersecting lattice paths and generating functions for tableaux,
{\it Mem. Amer. Math. Soc.} {\bf 115} (1995), no. 552, vi+109 pp.

\bibitem{kratt2}
C. Krattenthaler, A. J. Guttmann and X. G. Viennot, Vicious
walkers, friendly walkers and Young tableaux II: With a wall, {\it
J. Phys. A: Math. Gen.} {\bf 33} (2000), 8835 -- 8866.

\bibitem{macdonald}
I. G. Macdonald, {\it Symmetric Functions and Hall Polynomials}, Oxford University Press,
New York/London, 1979.

\bibitem{mills}
W. H. Mills, D. P. Robbins and H. Rumsey, Alternating sign matrices and descending plane partitions, {\it J. Combin. Theory Ser. A} {\bf 34} (1983), no. 3, 340 -- 359.

\bibitem{proctor}
R. A. Proctor, Bruhat lattices, plane partition generating functions, and minuscule representations, {\it Europ. J. Combin.} {\bf 5}, (1984), no. 4,  331 -- 350.

\bibitem{proctor2}
R. A. Proctor, New symmetric plane partition identities from invariant theory work of DeConcini and Procesi, {\it Europ. J. Combin.} {\bf 11}, (1990), no. 3, 289 -- 300.

\bibitem{slater}
L. J. Slater, {\it Generalized hypergeometric functions}, Cambridge University Press, Cambridge 1966.


\bibitem{stan}
R. P. Stanley, {\it Enumerative combinatorics}, vol. 2, Cambridge
University Press, Cambridge 1999.

\bibitem{stem2}
J. R. Stembridge, Hall-Littlewood functions, plane partitions and Rogers-Ramanujan identities, {\it Trans. Amer. Math. Soc.} {\bf 319}, (1990), no. 2, 469 -- 498.

\bibitem{stem}
J. R. Stembridge, The enumeration of totally symmetric plane partitions,
{\it Adv. Math.} {\bf 111} (1995), no. 2, 227 -- 243.

\bibitem{zeil}
D. Zeilberger, Proof of the refined alternating sign matrix conjecture,
{\it New York Journal of Mathematics} {\bf 2} (1996), 59 -- 68, electronic.

\end{thebibliography}
\end{document}